\documentclass[11pt,a4paper]{article}

\usepackage[T1]{fontenc}
\usepackage[utf8]{inputenc}
\usepackage{lmodern}
\usepackage{amsmath,amssymb,amsthm,mathtools,bm}
\usepackage{geometry}
\usepackage{microtype}
\usepackage{ragged2e}
\usepackage{hyperref}
\hypersetup{hidelinks,pdftitle={A Hyperbolic Microscopic Convexity Principle}}
\usepackage{enumitem}
\usepackage{titlesec}

\allowdisplaybreaks[3]
\numberwithin{equation}{section}

\newtheorem{theorem}{Theorem}[section]
\newtheorem{lemma}{Lemma}[section]

\newtheorem*{claim}{Claim}

\newcommand{\Om}{\Omega}
\newcommand{\Sn}{\mathbb{S}^{n}}
\newcommand{\Hn}{\mathbb{H}^{n+1}}
\newcommand{\R}{\mathbb{R}}
\newcommand{\Ssym}{\mathcal{S}^{n}}
\newcommand{\gb}{\bar g}
\newcommand{\bn}{\bar\nabla}
\newcommand{\be}{\bar e}
\newcommand{\Rc}{\bar R}
\newcommand{\eps}{\varepsilon}
\newcommand{\ph}{\varphi}
\newcommand{\sg}{\sigma}
\newcommand{\abs}[1]{\lvert #1\rvert}

\newcommand{\ip}[2]{\langle #1,#2\rangle}

\newcommand{\keywords}[1]{%
  \begingroup
  \renewcommand{\thefootnote}{}%
  \footnotetext{\textit{Key words and phrases. }#1}%
  \endgroup
}

\titleformat{\section}
  {\normalfont\Large\bfseries\raggedright}
  {\thesection}
  {1em}
  {}

\title{A Hyperbolic Microscopic Convexity Principle}
\author{Tong Liu, Lu Xu, and Yutian Zhou}
\date{}

\begin{document}

\maketitle

\begin{abstract}

We prove a constant rank theorem for the horospherical Hessian
$U[\ph]=\nabla^2\ph-\frac{\abs{\nabla\ph}^2}{2\ph}g
+\frac12(\ph-\ph^{-1})g$
of positive solutions to fully nonlinear elliptic equations on the round sphere.
The rescaled tensor $\varphi^{-1}U[\varphi]$ is Codazzi for the metric
$\varphi^{-2}g$.
Under a Bian--Guan type convexity condition with the hyperbolic position held fixed,
we show that if $U[\varphi]\geq 0$, then it is either positive definite everywhere or identically zero.
We apply this result to the $L_p$-Christoffel--Minkowski equation under an explicit convexity condition on the data.

\end{abstract}


\section{Introduction}

Constant rank theorems for convex solutions go back to Caffarelli and Friedman [1] and Korevaar and Lewis [2]. Caffarelli, Guan and Ma [5] developed a general method for fully nonlinear equations. Bian and Guan [3] considered equations in which the operator depends on the Hessian, the gradient, the value of the solution and the point \(x\). Codazzi tensors arise naturally in the hyperbolic Christoffel–Minkowski problem. Espinar, Gálvez and Mira [7] related the generalized hyperbolic Christoffel problem to conformally invariant equations on the sphere. The correspondence between horospherically convex hypersurfaces and conformal metrics and constant-rank problems for Codazzi curvature tensors have been studied extensively in the literature, we refer to  [8], [6]. On the hyperbolic side, Li and Xu [9] introduced the horospherical \(p\)-Brunn–Minkowski theory and the associated \(p\)-Christoffel–Minkowski problem. Luo and Wei [10] and Hu, Li and Xu [11] proved full-rank results for equations arising from this problem. Lei, Xu and Zhou proved a full-rank theorem for the hyperbolic \(k\)-Hessian equation \(\sigma_k(U[\varphi])=\varphi^{-k}f\) in [12] and later treated the Hessian quotient equation \(\sigma_k(U[\varphi])/\sigma_l(U[\varphi])=\varphi^{p-(k-l)}f\) in [13]. Both proofs work directly with \(U[\varphi]\) and use identities specific to the operator in question.

The main focus of this paper is on studying a general equation \(F(U,\nabla\varphi,\varphi,x)=0\). More specifically, we use a conformal change to let the rescaled tensor \(\widetilde U=\varphi^{-1}U[\varphi]\) satisfy the Codazzi identity. We use the hyperbolic position \(X\), which records both \(\varphi\) and \(\nabla\varphi\) at each point. Under the resulting convexity condition on the transformed operator, we prove that \(U[\varphi]\) has constant rank whenever it is nonnegative along a positive solution.

Let us introduce some notation. Let $\Om$ be a connected domain in $\Sn$ and denote by $\Ssym$ the space of real symmetric $n\times n$ matrices. Let $g$ be the round metric on $\Sn$,  and let $\ph>0$.  All derivatives appearing in
\begin{equation}\label{eq:U}
U[\ph]=\nabla^2\ph-\frac{\abs{\nabla\ph}^2}{2\ph}g
 +\frac12(\ph-\ph^{-1})g
\end{equation}
are taken with respect to $g$. Consider
\begin{equation}\label{eq:F}
F(U,\nabla\ph,\ph,x)=0.
\end{equation}
Put $\widetilde{U}=\ph^{-1}U[\ph]$. Section~2 defines the hyperbolic position $X$ and the transformed operator $\widetilde{F}$. Lemma~\ref{lem:equiv} proves
\begin{equation}\label{eq:equivintro}
F(U,\nabla\ph,\ph,x)=0\quad\Longleftrightarrow\quad \widetilde{F}(\widetilde{U},X,x)=0.
\end{equation}

\begin{theorem}\label{thm:main}
Let $n\ge2$ and suppose
\[
F=F(r,p,z,x)\in C^{2,1}\bigl(\Ssym\times\R^n\times(0,\infty)\times\Om\bigr),
\qquad F^{ij}=\left(\frac{\partial F}{\partial r_{ij}}\right)>0.
\]
Assume that
\begin{equation}\label{eq:struct}
\widetilde{F}(A^{-1},Y,x)\quad\text{is locally convex in $(A,x)$ for every fixed $Y$.}
\end{equation}
If $\ph\in C^{2,1}(\Om)$ is a positive solution of \eqref{eq:F} such that
\[
U[\ph]\ge0,
\]
then $U[\ph]$ has constant rank. Moreover,
\[
U[\ph]>0\quad\text{throughout $\Om$},\qquad\text{or}\qquad U[\ph]\equiv0.
\]
\end{theorem}

For comparison, the Euclidean theorem of Bian and Guan fixes the gradient variable $p$ and assumes the local convexity of $(A,z,x)\mapsto F(A^{-1},p,z,x)$. In \eqref{eq:struct}, we use the metric $g_Y$ defined in
Section~2.2. For each fixed $Y$, local convexity means that
$t\mapsto\widetilde{F}(A(t)^{-1},Y,x(t))$ is locally convex whenever
$x(t)$ is a $g_Y$-geodesic and $A(t)$ is positive definite with
entries affine in $t$ in a parallel $g_Y$-orthonormal frame
along $x(t)$. Once we fix $Y$, \eqref{eq:fixedY} gives both $z$ and $p$.
When we commute covariant derivatives, the resulting curvature term
has a favorable sign. We also need to keep track of how the position
variable changes and how the metric $\bar g$ associated with the
solution differs from the metric $g_Y$ obtained by fixing $Y$.

The rest of the paper is organized as follows. In Section~2,  we prove the conformal Codazzi identity and
rewrite the equation in terms of $\widetilde{U}=\ph^{-1}U[\ph]$.
We also derive the formulas for the position variable and the two
metrics that we use in Section~3 to prove Theorem~\ref{thm:main}.
In Section~4, we apply the theorem to the
$L_p$-Christoffel--Minkowski equation.

\section{Preliminary identities and analytic tools}
\subsection{The conformal Codazzi tensor}
Put
\begin{equation}\label{eq:Bg}
\gb=\ph^{-2}g,\qquad \widetilde{U}=\ph^{-1}U[\ph].
\end{equation}
Throughout the geometric part of the paper, $\nabla$ is the Levi-Civita connection of $g$, and $\bn$ is the Levi-Civita connection of $\gb$. The same bar convention is used for gradients and Hessians.

Let $\{e_i\}_{i=1}^n$ be a local $g$-orthonormal frame and set
\begin{equation}\label{eq:bar-frame}
  \bar e_i = \ph e_i.
\end{equation}
Then $\{\bar e_i\}_{i=1}^n$ is orthonormal with respect to $\bar g$, and
\begin{equation}\label{eq:Bgmatrix}
\widetilde{U}(\be_i,\be_j)=\ph U(e_i,e_j).
\end{equation}
Thus $U$ and $\widetilde{U}$ have the same eigendirections, nullspace and rank. When we differentiate $\widetilde U$ with respect to $\bar g$,
we choose a $\bar g$-orthonormal frame that is normal at the point
of calculation and diagonalizes $\widetilde U$ there.

\begin{lemma}\label{lem:codazzi}
The tensor $\widetilde{U}$ is Codazzi with respect to $\gb$. In a $\gb$-orthonormal frame diagonalizing $\widetilde{U}$, the sectional curvature is
\begin{equation}\label{eq:sectional}
K_{\gb}(\be_i,\be_j)=1+\lambda_i(\widetilde{U})+\lambda_j(\widetilde{U}).
\end{equation}
Consequently, if $\widetilde{U}\ge0$, then $K_{\gb}\ge1$.
\end{lemma}

\begin{proof}
Dividing \eqref{eq:U} by $\ph$ gives
\begin{equation}\label{eq:Bexplicit}
\widetilde{U}=\frac{\nabla^2\ph}{\ph}-\frac{\abs{\nabla\ph}^2}{2\ph^2}g
+\frac12(1-\ph^{-2})g.
\end{equation}
 Write
\begin{equation}\label{eq:Upsi}
  U_{ij} = \ph_{ij} + \ph g_{ij} - \psi g_{ij},
  \qquad
  \psi = \frac{\abs{\nabla\ph}^2}{2\ph}
       + \frac12(\ph + \ph^{-1}).
\end{equation}
Due to the Ricci identity, we have the following commutation formula on $\mathbb S^n$:
\[
\nabla_k\ph_{ij}-\nabla_j\ph_{ik}=R_{kji}{}^p\ph_p
=\ph_j g_{ik}-\ph_k g_{ij}.
\]
Then we have  $\nabla_k(\ph_{ij}+\ph g_{ij})-\nabla_j(\ph_{ik}+\ph g_{ik})=0$.  Differentiating $\psi$  gives
\begin{equation}\label{eq:psi-gradient}
\begin{aligned}
\psi_k&=\frac{\ph_p\ph_{pk}}{\ph}
-\frac{\abs{\nabla\ph}^2\ph_k}{2\ph^2}
+\frac12(1-\ph^{-2})\ph_k\\
&=\frac{\ph^p}{\ph}\left[\ph_{pk}-\frac{\abs{\nabla\ph}^2}{2\ph}g_{pk}
+\frac12(\ph-\ph^{-1})g_{pk}\right]
=\frac{\ph^p}{\ph}U_{pk}.
\end{aligned}
\end{equation}
Using \eqref{eq:psi-gradient}, direct computation yields 
\begin{equation}\label{eq:defect}
\nabla_kU_{ij}-\nabla_jU_{ik}
=g_{ik}\frac{\ph^p}{\ph}U_{pj}-g_{ij}\frac{\ph^p}{\ph}U_{pk}.
\end{equation}
Applying the Koszul formula to \(g\) and \(\bar g=\varphi^{-2}g\) gives the following identities.
\begin{align*}
2g(\nabla_XY,Z)={}&X\{g(Y,Z)\}+Y\{g(Z,X)\}-Z\{g(X,Y)\}\\
&-g(X,[Y,Z])+g(Y,[Z,X])+g(Z,[X,Y]),
\end{align*}
whereas the Koszul formula for $\gb=\ph^{-2}g$ is
\begin{align*}
2\gb(\bn_XY,Z)={}&X\{\gb(Y,Z)\}+Y\{\gb(Z,X)\}-Z\{\gb(X,Y)\}\\
&-\gb(X,[Y,Z])+\gb(Y,[Z,X])+\gb(Z,[X,Y]).
\end{align*}
The terms in the second identity are related to those in the first one by
\begin{align*}
X\{\gb(Y,Z)\}&=\ph^{-2}X\{g(Y,Z)\}-2\ph^{-3}X(\ph)g(Y,Z),\\
Y\{\gb(Z,X)\}&=\ph^{-2}Y\{g(Z,X)\}-2\ph^{-3}Y(\ph)g(Z,X),\\
Z\{\gb(X,Y)\}&=\ph^{-2}Z\{g(X,Y)\}-2\ph^{-3}Z(\ph)g(X,Y),\\
\gb(X,[Y,Z])&=\ph^{-2}g(X,[Y,Z]),
\end{align*}
and the last relation holds in exactly the same way for the other two bracket terms. Substitution into the $\gb$-formula shows directly that
\begin{align*}
2\gb(\bn_XY,Z)={}&2\ph^{-2}g(\nabla_XY,Z)\\
&-2\ph^{-3}\{X(\ph)g(Y,Z)+Y(\ph)g(X,Z)-Z(\ph)g(X,Y)\}.
\end{align*}
Since $\ph^{-2}g(\nabla_XY,Z)=\gb(\nabla_XY,Z)$, 
\begin{align*}
\gb(\bn_XY-\nabla_XY,Z)={}&-\frac{X(\ph)}{\ph}\gb(Y,Z)
-\frac{Y(\ph)}{\ph}\gb(X,Z)\\
&+\frac{Z(\ph)}{\ph}\gb(X,Y).
\end{align*}
Because this holds for every $Z$,
\begin{equation}\label{eq:connection}
\bn_XY=\nabla_XY-\frac{X(\ph)}{\ph}Y-\frac{Y(\ph)}{\ph}X
+\frac{g(X,Y)}{\ph}\nabla\ph.
\end{equation}
Here $\nabla\ph$ is the gradient with respect to $g$, which explains the factor $g(X,Y)$ in the last term. 
Then for vector fields $X,Y,Z$,
\begin{align*}
(\bn_X\widetilde{U})(Y,Z)&=X\{\widetilde{U}(Y,Z)\}-\widetilde{U}(\bn_XY,Z)-\widetilde{U}(Y,\bn_XZ),\\
(\nabla_X\widetilde{U})(Y,Z)&=X\{\widetilde{U}(Y,Z)\}-\widetilde{U}(\nabla_XY,Z)-\widetilde{U}(Y,\nabla_XZ).
\end{align*}
Subtracting these identities gives
\begin{equation}\label{eq:Bconn}
\begin{aligned}
(\bn_X\widetilde{U})(Y,Z)={}&(\nabla_X\widetilde{U})(Y,Z)\\
&-\widetilde{U}(\bn_XY-\nabla_XY,Z)-\widetilde{U}(Y,\bn_XZ-\nabla_XZ).
\end{aligned}
\end{equation}
At the point of calculation, take $X=e_k$, $Y=e_i$ and $Z=e_j$. Formula \eqref{eq:connection} gives
\[
\bn_{e_k}e_i-\nabla_{e_k}e_i
=-\frac{\ph_k}{\ph}e_i-\frac{\ph_i}{\ph}e_k+\frac{g_{ki}}{\ph}\nabla\ph.
\]
Since $\widetilde{U}=\ph^{-1}U$, the first term on the right of \eqref{eq:Bconn} is
\begin{equation}\label{eq:Bconn-gradient}
(\nabla_{e_k}\widetilde{U})(e_i,e_j)=\ph^{-1}\nabla_k U_{ij}-\ph^{-2}\ph_k U_{ij}.
\end{equation}
The correction from the first vector argument is
\begin{equation}\label{eq:Bconn-correction-one}
-\widetilde{U}(\bn_{e_k}e_i-\nabla_{e_k}e_i,e_j)
=\ph^{-2}(\ph_kU_{ij}+\ph_iU_{kj}-g_{ki}\ph^pU_{pj}),
\end{equation}
and the correction from the second vector argument is
\begin{equation}\label{eq:Bconn-correction-two}
-\widetilde{U}(e_i,\bn_{e_k}e_j-\nabla_{e_k}e_j)
=\ph^{-2}(\ph_kU_{ij}+\ph_jU_{ik}-g_{kj}\ph^pU_{ip}).
\end{equation}
Adding \eqref{eq:Bconn-gradient}, \eqref{eq:Bconn-correction-one}, and \eqref{eq:Bconn-correction-two} gives
\[
\bn_k\widetilde{U}_{ij}=\ph^{-1}\left\{\nabla_kU_{ij}
+\frac{\ph_k}{\ph}U_{ij}+\frac{\ph_i}{\ph}U_{kj}
+\frac{\ph_j}{\ph}U_{ik}
-\frac{\ph^p}{\ph}g_{ki}U_{pj}
-\frac{\ph^p}{\ph}g_{kj}U_{ip}\right\}.
\]
Consequently,
\[
\bn_k\widetilde{U}_{ij}-\bn_j\widetilde{U}_{ik}
=\ph^{-1}\left\{\nabla_kU_{ij}-\nabla_jU_{ik}
+g_{ij}\frac{\ph^p}{\ph}U_{pk}-g_{ik}\frac{\ph^p}{\ph}U_{pj}\right\}=0
\]
by \eqref{eq:defect}. Thus \(\widetilde U\) is Codazzi with respect to \(\bar g\).

It remains to compute the curvature of $\gb$ directly. We use the convention
\[
R(X,Y)Z=\nabla_X\nabla_YZ-\nabla_Y\nabla_XZ-\nabla_{[X,Y]}Z,
\qquad R_{ijkl}=g(R(e_i,e_j)e_l,e_k).
\]
For this calculation, write $\mathcal R_{ijkl}$ for the components
of the curvature tensor of $\bar g$ in the frame $e_i$. Then
\[
\mathcal R_{ijkl}
= \bar g\bigl(\bar R(e_i,e_j)e_l,e_k\bigr)
= \ph^{-2}g\bigl(\bar R(e_i,e_j)e_l,e_k\bigr).
\]
Let \(S\) be the difference between the two connections, so that
\[
S(X,Y)=\bn_XY-\nabla_XY
=-\frac{X(\ph)}{\ph}Y-\frac{Y(\ph)}{\ph}X+\frac{g(X,Y)}{\ph}\nabla\ph.
\]
 $\Rc(X,Y)Z=\bn_X\bn_YZ-\bn_Y\bn_XZ-\bn_{[X,Y]}Z$ and using $\bn_XY=\nabla_XY+S(X,Y)$ gives
\begin{align*}
\Rc(X,Y)Z={}&R(X,Y)Z+(\nabla_XS)(Y,Z)-(\nabla_YS)(X,Z)\\
&+S(X,S(Y,Z))-S(Y,S(X,Z)).
\end{align*}

At the point of calculation choose a $g$-normal frame. Then
\[
S(e_i,e_j)=-\frac{\ph_i}{\ph}e_j-\frac{\ph_j}{\ph}e_i
+\frac{\ph^p}{\ph}g_{ij}e_p,
\]
\[
\nabla_i\left(\frac{\ph_k}{\ph}\right)
=\frac{\ph_{ik}}{\ph}-\frac{\ph_i\ph_k}{\ph^2}.
\]
\[
\begin{aligned}
\mathcal R_{ijkl}
  &= \gb(\Rc(e_i,e_j)e_l,e_k)
   = \ph^{-2}g(\Rc(e_i,e_j)e_l,e_k) \\
  &= \ph^{-2}R_{ijkl}
     \begin{aligned}[t]
       &{}+ \ph^{-2}g\bigl(
         (\nabla_{e_i}S)(e_j,e_l)
         - (\nabla_{e_j}S)(e_i,e_l), e_k\bigr) \\
       &{}+ \ph^{-2}g\bigl(
         S(e_i,S(e_j,e_l))
         - S(e_j,S(e_i,e_l)), e_k\bigr).
     \end{aligned}
\end{aligned}
\]
Since \[
\nabla_i\left(\frac{\ph_j}{\ph}\right)=\nabla_j\left(\frac{\ph_i}{\ph}\right),
\] we have
\begin{equation}\label{eq:curvaturefour}
\begin{alignedat}{2}
\mathcal R_{ijkl} = \ph^{-2}\biggl\{ R_{ijkl}
&{}+ \left(\frac{\ph_{ik}}{\ph}
    - \frac{\abs{\nabla\ph}^2}{2\ph^2}g_{ik}\right)g_{jl}
&\quad &{}+ \left(\frac{\ph_{jl}}{\ph}
    - \frac{\abs{\nabla\ph}^2}{2\ph^2}g_{jl}\right)g_{ik} \\
&{}- \left(\frac{\ph_{il}}{\ph}
    - \frac{\abs{\nabla\ph}^2}{2\ph^2}g_{il}\right)g_{jk}
&\quad &{}- \left(\frac{\ph_{jk}}{\ph}
    - \frac{\abs{\nabla\ph}^2}{2\ph^2}g_{jk}\right)g_{il}
\biggr\}.
\end{alignedat}
\end{equation}
Choose the \(g\)-orthonormal frame so that \(\widetilde U\) is diagonal at the point. By (2.5), \((\varphi_{ij})\) is diagonal in the same frame. Since \(\bar e_i=\varphi e_i\), formula \eqref{eq:curvaturefour} gives
\[
K_{\gb}(\be_i,\be_j)=\ph^4\mathcal R_{ijij}
=\ph^2+\ph(\ph_{ii}+\ph_{jj})-\abs{\nabla\ph}^2.
\]
If $\lambda_i=\widetilde{U}(\be_i,\be_i)$, then \eqref{eq:Bexplicit} gives
\[
\lambda_i=\ph\ph_{ii}-\frac12\abs{\nabla\ph}^2+\frac12(\ph^2-1).
\]
Consequently,
\[
\ph^2+\ph(\ph_{ii}+\ph_{jj})-\abs{\nabla\ph}^2=1+\lambda_i+\lambda_j,
\]
which proves \eqref{eq:sectional}.
\end{proof}

\subsection{The transformed equation}
Write $\R^{n+1,1}=\R^{n+1}\times\R$ with
\[
\ip{(y,t)}{(y',t')}=y\cdot y'-tt',
\]
identify $\Sn$ with the unit sphere in $\R^{n+1}$, and put
\[
L(x)=(x,1),\qquad \Hn=\{Y\in\R^{n+1,1}:\ip{Y}{Y}=-1,\ Y^{n+2}>0\}.
\]
For fixed $Y\in\Hn$, define
\begin{equation}\label{eq:fixedY}
z=-\ip{Y}{L(x)},\qquad p=\nabla z,\qquad g_Y=z^{-2}g.
\end{equation}
If $e_i$ is a $g$-orthonormal frame, then $e_i^Y=ze_i$ is $g_Y$-orthonormal. Expressing the matrix $A$ with respect to the frame $e_i^Y$, define
\begin{equation}\label{eq:tildeF}
\widetilde{F}(A,Y,x)=F\bigl(z^{-1}A,p,z,x\bigr).
\end{equation}

The hyperbolic position associated with $\ph$ is given by
\begin{equation}\label{eq:X}
X=\frac12\ph(-x,1)
+\frac12\left(\frac{\abs{\nabla\ph}^2}{\ph}+\frac1\ph\right)(x,1)
-(\nabla\ph,0).
\end{equation}
We write $\bn_i X$ for its covariant
derivative along $\be_i$, using the Levi-Civita connection on $\Hn$.

\begin{lemma}\label{lem:position}
For every $x\in\mathbb{S}^n$, we have $X(x)\in\mathbb{H}^{n+1}$, and 
\begin{equation}\label{eq:Xrec}
\ph=-\ip{X}{L(x)},\qquad \ph_i=-\ip{X}{(e_i,0)}.
\end{equation}
Let
\[
E_j=-(e_j,0)+\frac{\ph_j}{\ph}L(x),
\qquad L(x)=(x,1).
\]
Using the frame $\be_i$ in \eqref{eq:bar-frame}, we have
\begin{equation}\label{eq:derX}
\bn_i X=\widetilde{U}_{ij}E_j.
\end{equation}
At a point where the $\gb$-normal frame diagonalizes $\widetilde{U}$, one has
\begin{equation}\label{eq:secX}
X_{;i}=\widetilde{U}_{ii}E_i,\qquad X_{;ii}=\widetilde{U}_{ii;j}E_j+\widetilde{U}_{ii}\bn_iE_i.
\end{equation}
Moreover,
\begin{equation}\label{eq:derE}
\bn_iE_i=(\nabla\ph,0)+\ph(x,0)
+\left(\ph_{ii}-\frac{\abs{\nabla\ph}^2}{\ph}\right)L(x)-\widetilde{U}_{ii}X.
\end{equation}
Here $x\in\Sn$, $L(x)=(x,1)$ is its null lift, $X$ is the position map defined in \eqref{eq:X}, $e_i\in T_x\Sn$ is $g$-orthonormal, and $(e_i,0)\in\R^{n+1,1}$ is the ambient lift of $e_i$. Moreover, $\ph_i=\nabla_{e_i}\ph$, $\be_i=\ph e_i$, $\bn_i=\bn_{\be_i}$, and $\widetilde{U}_{ij}=\widetilde{U}(\be_i,\be_j)$. In particular, $\bn_ZX=0$ whenever $Z\in\ker \widetilde{U}$.
\end{lemma}
\begin{proof}
In the following calculations, $\nabla_i$ applied to an $\R^{n+1,1}$-valued expression denotes ordinary differentiation in the direction $e_i$. The elementary Lorentzian products are
\[
\ip{(-x,1)}{(-x,1)}=0,\quad
\ip{(x,1)}{(x,1)}=0,\quad
\ip{(-x,1)}{(x,1)}=-2,
\]
and, because $\nabla\ph\perp x$,
\[
\ip{(\nabla\ph,0)}{(-x,1)}=\ip{(\nabla\ph,0)}{(x,1)}=0,\qquad
\ip{(\nabla\ph,0)}{(\nabla\ph,0)}=\abs{\nabla\ph}^2.
\]
Consequently
\[
\ip{X}{X}=-\ph\left(\frac{\abs{\nabla\ph}^2}{\ph}+\frac1\ph\right)
+\abs{\nabla\ph}^2=-\abs{\nabla\ph}^2-1+\abs{\nabla\ph}^2=-1.
\]
Moreover,
\[
X^{n+2}=\frac{\ph}{2}+\frac12\left(\frac{\abs{\nabla\ph}^2}{\ph}+\frac1\ph\right)>0,
\]
so $X\in\Hn$. Since $L(x)=(x,1)$,
\[
\ip{(-x,1)}{L(x)}=-2,\quad\ip{(x,1)}{L(x)}=0,\quad
\ip{(\nabla\ph,0)}{L(x)}=0,
\]
which gives
\[
-\ip{X}{L(x)}=-\frac\ph2\ip{(-x,1)}{(x,1)}=\ph.
\]
For each unit tangent vector $e_i$, 
\[
-\ip{X}{(e_i,0)}
=\ip{(\nabla\ph,0)}{(e_i,0)}
=\ph_i.
\]
This proves \eqref{eq:Xrec}.

We next differentiate $X$. Choose a $g$-orthonormal frame
$\{e_i\}$ that is normal at the point of calculation. At this point,
\[
\nabla_i(x,1)=(e_i,0),\qquad
\nabla_i(-x,1)=-(e_i,0),\qquad
\nabla_i(\nabla\ph,0)=(\ph_{ij}e_j-\ph_i x,0).
\]
Differentiating \eqref{eq:X} and using $(-x,1)=L(x)-2(x,0)$,
we obtain
\begin{align*}
\nabla_i X
&=-U_{ij}(e_j,0)
 +\left(
    \frac{\ph_i}{2}
    +\frac{\ph_j\ph_{ji}}{\ph}
    -\frac{\abs{\nabla\ph}^2\ph_i}{2\ph^2}
    -\frac{\ph_i}{2\ph^2}
  \right)L(x)\\
&=-U_{ij}(e_j,0)+U_{ij}\frac{\ph_j}{\ph}L(x)=U_{ij}E_j,
\end{align*}
where
\[
E_j=-(e_j,0)+\frac{\ph_j}{\ph}L(x).
\]
Since $\be_i=\ph e_i$, it follows that
$\bn_iX=\ph U_{ij}E_j$. By \eqref{eq:Bgmatrix},
$\ph U_{ij}=\widetilde U(\be_i,\be_j)$, which gives
\eqref{eq:derX}.
For the next calculation, keep the values of $\{\be_i\}$ at this
point and extend them to a local $\gb$-orthonormal frame that is
normal there. Differentiating \eqref{eq:derX} gives
\[
X_{;ii}=\widetilde{U}_{ij;i}E_j+\widetilde{U}_{ij}\bn_iE_j.
\]
Since $\widetilde{U}$ is Codazzi and diagonal at the point,
\[
\widetilde{U}_{ij;i}=\widetilde{U}_{ii;j},\qquad \widetilde{U}_{ij}\bn_iE_j=\widetilde{U}_{ii}\bn_iE_i.
\]
This proves \eqref{eq:secX}.

Direct differentiation of the formula for $E_i$ gives
\begin{equation}\label{eq:ambient-Ei-derivative}
dE_i(\be_i)=(\nabla\ph,0)+\ph(x,0)
+\left(\ph_{ii}-\frac{\abs{\nabla\ph}^2}{\ph}\right)L(x).
\end{equation}
Since $\ip{X}{X}=-1$, the orthogonal projection of an ambient vector $Z$ onto $T_X\Hn$ is $Z+\ip{Z}{X}X$. Using $X_{;i}=\widetilde{U}_{ii}E_i$ and $\ip{E_i}{E_i}=1$, we obtain
\begin{flalign*}
&\hspace{12em}
\begin{aligned}
\bn_i E_i
  &= dE_i(\be_i) + \ip{dE_i(\be_i)}{X}X \\
  &= dE_i(\be_i) - \ip{E_i}{X_{;i}}X
   = dE_i(\be_i) - \widetilde{U}_{ii}X,
\end{aligned}
&&
\end{flalign*}
which, together with \eqref{eq:ambient-Ei-derivative}, gives \eqref{eq:derE}.
\end{proof}
\begin{lemma}\label{lem:equiv}
Along the solution, the original equation \eqref{eq:F} is equivalent to
\begin{equation}\label{eq:tildeFequiv}
\widetilde{F}(\widetilde{U},X,x)=0,
\end{equation}
where $X$ is the position map defined in \eqref{eq:X}, and the matrix argument $\widetilde{U}$ is measured in a $\gb$-orthonormal frame.
\end{lemma}
\begin{proof}
Fix $x$ and take a $g$-orthonormal frame $e_i$. Use the associated $\gb$-orthonormal frame $\be_i=\ph e_i$, and denote the two matrix representations by $[U]_g$ and $[\widetilde{U}]_{\gb}$. By \eqref{eq:Bgmatrix},
\[
[\widetilde{U}]_{\gb}=\ph[U]_g,\qquad [U]_g=\ph^{-1}[\widetilde{U}]_{\gb}.
\]
\begin{align*}
\widetilde{F}([\widetilde{U}]_{\gb},X,x)
&=F\Bigl(\frac{[\widetilde{U}]_{\gb}}{-\ip{X}{L(x)}},
-\sum_i\ip{X}{(e_i,0)}e_i,-\ip{X}{L(x)},x\Bigr)\\
&=F\bigl(\ph^{-1}[\widetilde{U}]_{\gb},\nabla\ph,\ph,x\bigr)\\
&=F\bigl([U]_g,\nabla\ph,\ph,x\bigr)\\
&=F(U[\ph],\nabla\ph,\ph,x)=0,
\end{align*}
where the first equality is \eqref{eq:tildeF}, the second uses \eqref{eq:Xrec}, and the third uses $\ph^{-1}[\widetilde{U}]_{\gb}=[U]_g$.
\end{proof}

\begin{lemma}\label{lem:freeze}
Fix $x_0\in\Om$ and freeze the position variable at $Y=X(x_0)\in\Hn$. Put
\[
z(x)=-\ip{X(x_0)}{L(x)},\qquad g_Y=z^{-2}g,
\]
and denote the Levi-Civita connection of $g_Y$ by $D$. At $x_0$,
\begin{equation}\label{eq:equalconn}
z=\ph,\quad dz=d\ph,\quad g_Y=\gb,\quad D=\bn.
\end{equation}
If $\be_i$ is an eigenvector of $\widetilde{U}$ at $x_0$, then
\begin{equation}\label{eq:connestimate}
\
\left|\bar\nabla_{\bar e_i}(\bar\nabla-D)\right|
\le C|\widetilde U_{ii}|.
\
\end{equation}
The derivative and norm are evaluated at $x_0$, and $C$ is locally uniform.
\end{lemma}
\begin{proof}
By \eqref{eq:Xrec}, at $x_0$,
\[
z=-\ip{X(x_0)}{L(x_0)}=\ph,\qquad z_i=-\ip{X(x_0)}{(e_i,0)}=\ph_i.
\]
Consequently
\[
(g_Y)_{ij}=z^{-2}g_{ij}=\ph^{-2}g_{ij}=\gb_{ij},
\qquad \nabla_k(g_Y)_{ij}=-2z^{-3}z_kg_{ij}=-2\ph^{-3}\ph_kg_{ij}=\nabla_k\gb_{ij}.
\]
The Koszul formula gives $D=\bn$ at $x_0$, and proves \eqref{eq:equalconn}.
We use $\nabla_i\nabla_j x=-g_{ij}x$
and evaluate at $x_0$ to obtain
\begin{equation}\label{eq:fixed-z-hessian}
z_{ij}=\ip{X(x_0)}{(x,0)}g_{ij}
=\left\{\frac{\abs{\nabla\ph}^2}{2\ph}-\frac12(\ph-\ph^{-1})\right\}g_{ij}.
\end{equation}
Solving \eqref{eq:U} for the Hessian gives
\begin{equation}\label{eq:fixed-phi-hessian}
\ph_{ij}=U_{ij}+\left\{\frac{\abs{\nabla\ph}^2}{2\ph}
-\frac12(\ph-\ph^{-1})\right\}g_{ij}.
\end{equation}
Subtracting \eqref{eq:fixed-z-hessian} from \eqref{eq:fixed-phi-hessian} gives
\begin{equation}\label{eq:fixed-hessian-difference}
\ph_{ik}-z_{ik}=U_{ik}.
\end{equation}
Since $z=\ph$ and $z_k=\ph_k$ at $x_0$, \eqref{eq:fixed-hessian-difference} yields
\begin{equation}\label{eq:freeze-component-general}
\nabla_i\left(\frac{\ph_k}{\ph}-\frac{z_k}{z}\right)
=\frac{\ph_{ik}-z_{ik}}{\ph}=\frac{U_{ik}}{\ph}
=\frac{\widetilde{U}(\be_i,\be_k)}{\ph^2}.
\end{equation}
At a point where $\widetilde{U}$ is diagonal, \eqref{eq:freeze-component-general} becomes
\begin{equation}\label{eq:freeze-component-diagonal}
\nabla_i\left(\frac{\ph_k}{\ph}-\frac{z_k}{z}\right)
=\begin{cases}0,&k\ne i,\\ \widetilde{U}_{ii}/\ph^2,&k=i.\end{cases}
\end{equation}
Subtracting formula \eqref{eq:connection} for $\gb=\ph^{-2}g$ and $g_Y=z^{-2}g$ gives, for tangent vectors $V_1,V_2$,
\begin{equation}\label{eq:connection-difference}
\begin{aligned}
(\bn-D)_{V_1}V_2={}&-\left(\frac{V_1(\ph)}{\ph}-\frac{V_1(z)}{z}\right)V_2
-\left(\frac{V_2(\ph)}{\ph}-\frac{V_2(z)}{z}\right)V_1\\
&+g(V_1,V_2)\left(\frac{\nabla\ph}{\ph}-\frac{\nabla z}{z}\right).
\end{aligned}
\end{equation}
At $x_0$, the right-hand side of \eqref{eq:connection-difference} vanishes by $z=\ph$ and $dz=d\ph$. Extend $V_1,V_2$ so that their covariant derivatives vanish at $x_0$. Differentiating \eqref{eq:connection-difference} in the $\be_i$ direction and using \eqref{eq:freeze-component-diagonal} gives
\[
\bigl\lvert\bn_{\be_i}(\bn-D)_{V_1}V_2\bigr\rvert\le C|\widetilde{U}_{ii}|\abs{V_1}\abs{V_2},
\]
which proves \eqref{eq:connestimate}.
\end{proof}

\subsection{Analytic tools for the constant-rank calculation}
We use a semicolon to denote covariant derivatives with respect to
the metric in use. In Section~3, we use the metric $\gb$.
At the point where we carry out the calculation, we choose a normal
frame in which $\widetilde{U}$ is diagonal.
\begin{lemma}\label{lem:quotient}
Let $A(x)\ge0$ be a smooth symmetric two-tensor and let \(l(l<n)\) be its minimum rank. Put
\[
A^\eps=A+\eps I,\qquad q_\eps(A)=\frac{\sg_{l+2}(A^\eps)}{\sg_{l+1}(A^\eps)},
\qquad \phi_\eps(A)=\sg_{l+1}(A^\eps)+q_\eps(A).
\]
Near a point of rank $l$, let $G$ contain the indices of the $l$ eigenvalues bounded below by a positive constant and let $B$ contain the remaining indices. Choose an orthonormal frame in which $A$ is diagonal. $\sg_k(G)$ is the $k$th elementary symmetric function of $\{A_{aa}+\eps:a\in G\}$, while $\sg_k(B\mid i)$ is the same function of $\{A_{rr}+\eps:r\in B,\ r\ne i\}$. Then there are constants $0<c\le C$, independent of $\eps$, such that
\begin{equation}\label{eq:phi-comparison}
c\sum_{r\in B}(A_{rr}+\eps)\le\phi_\eps(A)
\le C\sum_{r\in B}(A_{rr}+\eps).
\end{equation}
For $\alpha,\beta\in\{1,\ldots,n\}$,
\begin{equation}\label{eq:phi-second}
\begin{aligned}
&\sum_{i,j,k,m} q_\eps^{ij,km}
 A^\eps_{ij;\alpha} A^\eps_{km;\beta}
 = O\left(\phi_\eps+\sum_{i,j\in B}\abs{\nabla A^\eps_{ij}}\right)
 -\frac{1}{\sg_1(B)}
 \sum_{\substack{i,j\in B\\ i\ne j}}
 A^\eps_{ij;\alpha}A^\eps_{ji;\beta} \\
&\quad
 -\frac{1}{\sg_1(B)^3}
 \sum_{i\in B}V_{i\alpha}V_{i\beta}
 -2\sum_{i\in B,j\in G}
 \frac{\sg_1(B\mid i)^2-\sg_2(B\mid i)}
      {\sg_1(B)^2(A_{jj}+\eps)}
 A^\eps_{ij;\alpha}A^\eps_{ij;\beta},
\end{aligned}
\end{equation}
where
\[
V_{i\alpha}
=\sg_1(B)A^\eps_{ii;\alpha}
 -(A_{ii}+\eps)\sum_{j\in B}A^\eps_{jj;\alpha}.
\]
Moreover, $\phi_\eps(A)$ converges locally in $C^1$
to a $C^{1,1}$ function.
\end{lemma}
\begin{proof}
The comparison in (2.32) follows from the definition of $\phi_\eps$
and the bounds on the eigenvalues indexed by $G$.
Formula (2.33) follows from the calculation in
[3, Proposition~2.1, equation~(2.7)], applied in a normal
orthonormal frame at a point where $A$ is diagonal.

For the regularity assertion, we use the proof of
[3, Corollary~2.2]. Since $A\ge 0$, the scalar function
$A(v,v)$ is nonnegative for every local vector field $v$.
Applying the gradient estimate for nonnegative $C^2$ functions
to $A(v,v)$ with $v=e_i+e_j$ and $v=e_i-e_j$, and using
polarization, gives
\[
  \lvert\nabla A_{ij}\rvert
  \le C\bigl(\sqrt{A_{ii}}+\sqrt{A_{jj}}\bigr)
\]
at the chosen point. With this estimate, the proof of
[3, Corollary~2.2] gives a local $C^{1,1}$ bound for $\phi_\eps$
that is uniform in $\eps$. Since $\phi_\eps$ converges pointwise,
it converges locally in $C^1$ to a $C^{1,1}$ function.
\end{proof}

\section{Proof of the main theorem}
The following differential inequality is the hyperbolic counterpart of Bian--Guan's Theorem~3.2 \cite{BG09}. 

\begin{theorem}\label{thm:differential}
Assume the hypotheses of Theorem~\ref{thm:main}, and let $x_0$ be a point where $\widetilde{U}=\ph^{-1}U[\ph]$ has minimum rank $l$. Let $\mathcal O$ be a small open neighborhood such that for each $x\in\mathcal O$, there are $l$ "good" eigenvalues of $\widetilde{U}$ which are bounded below by a positive constant, and the other $n-l$"bad" eigenvalues of $\widetilde{U}$ are very small. Put
\[
\phi=\sg_{l+1}(\widetilde{U})+
\begin{cases}
\sg_{l+2}(\widetilde{U})/\sg_{l+1}(\widetilde{U}),&\sg_{l+1}(\widetilde{U})>0,\\
0,&\sg_{l+1}(\widetilde{U})=0.
\end{cases}
\]
Then $\phi$ is locally $C^{1,1}$, and there is a constant $C$ such that
\begin{equation}\label{eq:main-differential}
\widetilde{F}^{\alpha\beta}\phi_{;\alpha\beta}
\le C\bigl(\phi+\abs{\bn\phi}\bigr).
\end{equation}
Consequently $\widetilde{U}$ has constant rank in $\mathcal O$. Moreover, on the connected set $\Om$, either $\widetilde{U}>0$ or $\widetilde{U}\equiv0$.
\end{theorem}
\begin{proof}
We now apply the preceding identities directly to the tensor $\widetilde{U}=\ph^{-1}U[\ph]$ and the position map $X$ occurring in the theorem. All covariant derivatives below are taken with respect to $\gb$ and are denoted by $\bn$; semicolons have the same meaning. Lemma~\ref{lem:codazzi} gives that $\widetilde{U}\ge0$ is Codazzi and that every sectional curvature of $\gb$ is at least one. Section~2.2 gives
\[
\widetilde{F}(\widetilde{U},X,x)=0.
\]
 If $l=n$, constant rank is already proved, so we only consider $l<n$. Denote $G$ the index set of these "good" eigenvalues and B the index set of "bad" eigenvalues. Put
\begin{equation}\label{eq:reg}
\widetilde{U}^\eps=\widetilde{U}+\eps\gb,\qquad
\phi_\eps=\sg_{l+1}(\widetilde{U}^\eps)+\frac{\sg_{l+2}(\widetilde{U}^\eps)}{\sg_{l+1}(\widetilde{U}^\eps)}.
\end{equation}
Lemma~\ref{lem:quotient}, applied to $\widetilde{U}^\eps$, gives on a smaller coordinate ball
\begin{equation}\label{eq:boundsreg}
\phi_\eps\ge0,\quad\phi_\eps(x_0)=O(\eps),\quad
c\sum_{i\in B}(\widetilde{U}_{ii}+\eps)\le\phi_\eps
\le C\sum_{i\in B}(\widetilde{U}_{ii}+\eps).
\end{equation}
In particular,
\begin{equation}\label{eq:bad-small}
\eps\le C\phi_\eps,\qquad \widetilde{U}_{ii}+\eps\le C\phi_\eps\quad(i\in B).
\end{equation}
These are the counterparts of \cite[(3.6)--(3.10)]{BG09}. Unlike the Euclidean regularization $D^2u+\eps I$, the metric term here is $\eps\gb$. The tensor $\widetilde{U}^\eps$ satisfies the perturbed equation
\begin{equation}\label{eq:residual}
R_\eps=\widetilde{F}(\widetilde{U}^\eps,X,x)-\widetilde{F}(\widetilde{U},X,x),\qquad \widetilde{F}(\widetilde{U}^\eps,X,x)=R_\eps.
\end{equation}
Since $\bn(\widetilde{U}+\eps\gb)=\bn \widetilde{U}$, differentiating this identity up to twice and using \eqref{eq:bad-small} gives
\begin{equation}\label{eq:residual-bd}
\abs{R_\eps}+\abs{\bn R_\eps}+\abs{\bn^2R_\eps}\le C\eps\le C\phi_\eps.
\end{equation}
Fix an arbitrary point $x_1\in\mathcal O$. We will prove \eqref{eq:main-differential} at $x_1$. By choosing a suitable $\gb$-normal frame, we may assume $\widetilde{U}(x_1)$ is diagonal. We write the regularized eigenvalues directly as $\widetilde{U}_{ii}+\eps$. For every derivative direction $\be_\alpha$, the notation used below is
\begin{equation}\label{eq:semi-def}
\widetilde{U}_{ij;\alpha}=(\bn_{\be_\alpha}\widetilde{U})(\be_i,\be_j),\qquad
\widetilde{U}_{ij;\alpha\beta}=(\bn_{\be_\beta}\bn_{\be_\alpha}\widetilde{U})(\be_i,\be_j).
\end{equation}
Since $\bn\gb=0$, the matrix chain rule gives
\begin{equation}\label{eq:chain1}
(\phi_\eps)_{;\alpha}=\phi^{ij}\widetilde{U}_{ij;\alpha},
\end{equation}
\begin{equation}\label{eq:chain2}
(\phi_\eps)_{;\alpha\beta}
=\phi^{ij}\widetilde{U}_{ij;\alpha\beta}+\phi^{ij,km}\widetilde{U}_{ij;\alpha}\widetilde{U}_{km;\beta}.
\end{equation}
Consequently,
\begin{equation}\label{eq:chain-contraction}
\widetilde{F}^{\alpha\beta}(\phi_\eps)_{;\alpha\beta}
=\phi^{ij}\widetilde{F}^{\alpha\beta}\widetilde{U}_{ij;\alpha\beta}
+\widetilde{F}^{\alpha\beta}\phi^{ij,km}\widetilde{U}_{ij;\alpha}\widetilde{U}_{km;\beta}.
\end{equation}

We first treat the fourth-order term in \eqref{eq:chain-contraction}. This is the point at which the spherical background differs from the Euclidean calculation. At the diagonal point the spectral derivative $\phi^{ij}$ vanishes for $i\ne j$, so this term is
\[
\sum_i\phi^{ii}\widetilde{F}^{\alpha\beta}\widetilde{U}_{ii;\alpha\beta}.
\]
Two Codazzi exchanges give
\[
\widetilde{U}_{ii;\alpha\beta}=\widetilde{U}_{i\alpha;i\beta},\qquad
\widetilde{U}_{\alpha\beta;ii}=\widetilde{U}_{\alpha i;\beta i}.
\]
For the curvature convention used in this paper,
\[
(\bn_k\bn_j-\bn_j\bn_k)\omega_i=\Rc_{kji}{}^p\omega_p.
\]
Applying this identity to the two covariant indices of $\widetilde{U}$ yields
\begin{equation}\label{eq:commutator}
\widetilde{U}_{ii;\alpha\beta}-\widetilde{U}_{\alpha\beta;ii}
=\Rc_{\beta ii}{}^{r}\widetilde{U}_{r\alpha}+\Rc_{\beta i\alpha}{}^{r}\widetilde{U}_{ir}.
\end{equation}
At the point under consideration $\widetilde{U}$ is diagonal. Formula \eqref{eq:curvaturefour} gives
\[
\Rc_{i\alpha i\beta}=\bar g\big(
\bar R(\bar e_i,\bar e_\alpha)\bar e_\beta,\bar e_i
\big)
=(\widetilde{U}_{ii}+\tfrac12)\delta_{\alpha\beta}
+(\widetilde{U}_{\alpha\beta}+\tfrac12\delta_{\alpha\beta})
-(\widetilde{U}_{i\beta}+\tfrac12\delta_{i\beta})\delta_{\alpha i}
-(\widetilde{U}_{\alpha i}+\tfrac12\delta_{\alpha i})\delta_{i\beta}.
\]
Thus $\Rc_{i\alpha i\beta}=0$ for $\alpha\ne\beta$, and contraction of \eqref{eq:commutator} gives
\begin{equation}\label{eq:comm-contract}
\widetilde{F}^{\alpha\beta}\widetilde{U}_{ii;\alpha\beta}
=\widetilde{F}^{\alpha\beta}\widetilde{U}_{\alpha\beta;ii}
-\sum_{a\ne i}\widetilde{F}^{aa}\Rc_{aiai}(\widetilde{U}_{aa}-\widetilde{U}_{ii}).
\end{equation}
Consequently the fourth-order part of \eqref{eq:chain-contraction} becomes
\begin{equation}\label{eq:fourth}
\sum_i\phi^{ii}\widetilde{F}^{\alpha\beta}\widetilde{U}_{ii;\alpha\beta} = \sum_i\phi^{ii}\widetilde{F}^{\alpha\beta}\widetilde{U}_{\alpha\beta;ii} - \sum_i\phi^{ii}\sum_{a\ne i}\widetilde{F}^{aa}\Rc_{aiai}(\widetilde{U}_{aa}-\widetilde{U}_{ii}).
\end{equation}
 Differentiate $\widetilde{F}(\widetilde{U}^\eps,X,x)=R_\eps$ in the direction $\be_i$ and use \eqref{eq:residual-bd} to obtain
\begin{equation}\label{eq:diff1}
\widetilde{F}^{ab}\widetilde{U}_{ab;i}+\widetilde{F}^Y\cdot X_{;i}+\widetilde{F}^{x_i}=O(\phi_\eps).
\end{equation}
Differentiating once more gives
\begin{equation}\label{eq:diff2}
\begin{aligned}
&\widetilde{F}^{ab}\widetilde{U}_{ab;ii}+\widetilde{F}^{ab,cd}\widetilde{U}_{ab;i}\widetilde{U}_{cd;i}
+2\widetilde{F}^{ab,Y}(\widetilde{U}_{ab;i},X_{;i})\\
&\quad+2\widetilde{F}^{ab,x_i}\widetilde{U}_{ab;i}+\widetilde{F}^{Y,Y}(X_{;i},X_{;i})
+2\widetilde{F}^{Y,x_i}\cdot X_{;i}\\
&\quad+\widetilde{F}^Y\cdot X_{;ii}+\widetilde{F}^{x_i,x_i}=O(\phi_\eps).
\end{aligned}
\end{equation}
 Solving \eqref{eq:diff2} for $\widetilde{F}^{ab}\widetilde{U}_{ab;ii}$, substituting it into \eqref{eq:fourth}, and then substituting the result into \eqref{eq:chain-contraction} gives 
\begin{equation}\label{eq:big-expansion}
\begin{aligned}
&\widetilde{F}^{\alpha\beta}(\phi_\eps)_{;\alpha\beta}
 = \sum\widetilde{F}^{\alpha\beta}\phi^{ij,km}
   \widetilde{U}_{ij;\alpha}\widetilde{U}_{km;\beta} \\
&\quad -\sum_i\phi^{ii}\Bigl\{
  2\widetilde{F}^{ab,Y}(\widetilde{U}_{ab;i},X_{;i})
  +\widetilde{F}^{Y,Y}(X_{;i},X_{;i})
  +2\widetilde{F}^{Y,x_i}\cdot X_{;i}
  +\widetilde{F}^Y\cdot X_{;ii}
 \Bigr\} \\
&\quad -\sum_i\phi^{ii}\sum_{a\ne i}
  \widetilde{F}^{aa}\Rc_{aiai}
  (\widetilde{U}_{aa}-\widetilde{U}_{ii})
  +O(\phi_\eps) \\
&\quad -\sum_i\phi^{ii}\Bigl\{
  \widetilde{F}^{ab,cd}\widetilde{U}_{ab;i}\widetilde{U}_{cd;i}
  +2\widetilde{F}^{ab,x_i}\widetilde{U}_{ab;i}
  +\widetilde{F}^{x_i,x_i}
 \Bigr\}.
\end{aligned}
\end{equation}
At the diagonal point $x_1$,
\begin{equation}\label{eq:phi-coeff}
\begin{aligned}
\phi^{ii}&=\sg_l(G)+\frac{\sg_1(B\mid i)^2-\sg_2(B\mid i)}{
\bigl(\sum_{r\in B}(\widetilde{U}_{rr}+\eps)\bigr)^2}+O(\phi_\eps)
&& (i\in B),\\
\phi^{jj}&=O(\phi_\eps)&& (j\in G).
\end{aligned}
\end{equation}
Formula \eqref{eq:phi-coeff} follows from \cite[Lemmas~2.3--2.4]{BG09}, applied to $\widetilde{U}^\eps$. Lemma~\ref{lem:position} gives, for every $i=1,\ldots,n$,
\begin{equation}\label{eq:X-in-proof}
\begin{aligned}
X_{;i}&=\widetilde{U}_{ij}E_j=\widetilde{U}_{ii}E_i,\\
X_{;ii}&=\bn_i(\widetilde{U}_{ij}E_j)=\widetilde{U}_{ij;i}E_j+\widetilde{U}_{ij}\bn_iE_j\\
&=\widetilde{U}_{ii;j}E_j+\widetilde{U}_{ii}\bn_iE_i.
\end{aligned}
\end{equation}
The last equality uses that $\widetilde{U}$ is diagonal at the point and $\widetilde{U}_{ij;i}=\widetilde{U}_{ii;j}$ by the Codazzi identity. Here
\[
E_j=-(e_j,0)+\frac{\ph_j}{\ph}L(x),\qquad L(x)=(x,1).
\]
Formula \eqref{eq:derE} gives $\abs{\bn_iE_i}\le C$. For $i\in B$, \eqref{eq:bad-small} and \eqref{eq:X-in-proof} give
\[
\abs{X_{;i}}\le C\phi_\eps,\qquad
\abs{\widetilde{U}_{ii}\bn_iE_i}\le C\phi_\eps\quad(i\in B).
\]
Consequently,
\begin{equation}\label{eq:X-error}
\begin{aligned}
&\sum_i\phi^{ii}\bigl\{2\widetilde{F}^{ab,Y}(\widetilde{U}_{ab;i},X_{;i})+\widetilde{F}^{Y,Y}(X_{;i},X_{;i})\\
&\qquad+2\widetilde{F}^{Y,x_i}\cdot X_{;i}+\widetilde{F}^Y\cdot X_{;ii}\bigr\}
=O(\phi_\eps+\sum_{i,j\in B}\abs{\bn \widetilde{U}_{ij}}).
\end{aligned}
\end{equation}
Indeed, for $i\in G$, $\phi^{ii}=O(\phi_\eps)$, so all corresponding terms are $O(\phi_\eps)$. For $i\in B$, every term containing $X_{;i}$ is $O(\phi_\eps)$. The only remaining term is
\[\sum_i
\phi^{ii}\widetilde{F}^Y\cdot \widetilde{U}_{ii;j}E_j=O(\phi_\eps)+\sum_{i\in B}\phi^{ii}\widetilde{F}^Y\cdot \widetilde{U}_{ii;j}E_j=O(\phi_\eps+\sum_{i,j\in B}\abs{\bn \widetilde{U}_{ij}}).
\]
At $x_1$, fix $Y=X(x_1)$, put
\[
z(x)=-\ip{Y}{L(x)},\qquad g_Y=z^{-2}g,
\]
and let $D$ be the Levi-Civita connection of $g_Y$. Equation \eqref{eq:equalconn} gives $D=\bn$ at this point. The part of the last line with $a,b,c,d\in G$ is therefore decomposed as
\begin{equation}\label{eq:decomposition}
\begin{aligned}
&\sum_{a,b,c,d\in G}\widetilde{F}^{ab,cd}\widetilde{U}_{ab;i}\widetilde{U}_{cd;i}
+2\sum_{a,b\in G}\widetilde{F}^{ab,x_i}\widetilde{U}_{ab;i}+\widetilde{F}^{x_i,x_i}\\
&=\sum_{a,b,c,d\in G}\widetilde{F}^{ab,cd}\widetilde{U}_{ab;i}\widetilde{U}_{cd;i}
+2\sum_{a,b\in G}\widetilde{F}_D^{ab,x_i}\widetilde{U}_{ab;i}+\widetilde{F}_D^{x_i,x_i}\\
&\quad+2\sum_{a,b\in G}(\widetilde{F}^{ab,x_i}-\widetilde{F}_D^{ab,x_i})\widetilde{U}_{ab;i}
+\widetilde{F}^{x_i,x_i}-\widetilde{F}_D^{x_i,x_i}.
\end{aligned}
\end{equation}
The first three terms on the right-hand side are calculated with $Y$ fixed and with the connection $D$, so condition \eqref{eq:struct} applies to them. The last two terms compare the actual connection $\bn$ with $D$; they are estimated in the following claim.

\begin{claim}
At the arbitrary point $x_1$ fixed above,
\begin{equation}\label{eq:claim}
\widetilde{F}^{ab,x_i}=\widetilde{F}_D^{ab,x_i},\qquad
\abs{\widetilde{F}^{x_i,x_i}-\widetilde{F}_D^{x_i,x_i}}
\le C\abs{\bn_{\be_i}(\bn-D)}\le C\widetilde{U}_{ii}.
\end{equation}
\end{claim}
\begin{proof}[Proof of the claim]
Keep $Y=X(x_1)$ fixed. At $x_1$, Lemma~\ref{lem:freeze} gives
$z=\ph$, $dz=d\ph$, and $D=\bn$. Choose a $\gb$-geodesic
$x(t)$ with $x(0)=x_1$ and $\dot x(0)=\be_i$, together with a
$\bn$-parallel $\gb$-orthonormal frame $E_a(t)$ satisfying
$E_a(0)=\be_a$. Since $D=\bn$ at $x_1$, $D_{\dot x}\dot x(0)=0$.

Let $A=\widetilde U^\eps(x_1)$, and denote its parallel transports
as a covariant two-tensor under $\bn$ and $D$ by $A_{\bn}(t)$ and
$A_D(t)$, respectively. Let
$(A_*)_{ab}(t)=A_*(t)(E_a(t),E_b(t))$, with $*\in\{\bn,D\}$.
Since $\bn=D$ at $x_1$ and $E_a(t)$ is $\bn$-parallel,
\[
(A_{\bn})_{ab}(0)=(A_D)_{ab}(0)=A_{ab},\qquad
(A_{\bn})'_{ab}(0)=(A_D)'_{ab}(0)=0.
\]

The matrix argument of $\widetilde F$ is measured in a
$g_Y$-orthonormal frame. Because $E_a(t)$ is $\gb$-orthonormal,
$(z/\ph)E_a(t)$ is $g_Y$-orthonormal. Thus we must evaluate
$\widetilde F$ on the matrix with entries
$(\bar A_*)_{ab}(t)=\bigl(\frac{z(x(t))^2}{\ph(x(t))^2}(A_*)_{ab}(t)\bigr)$.
Using $z=\ph$ and $dz=d\ph$ at $x_1$, we obtain
\[
\left.\frac{d}{dt}\left(\frac{z^2}{\ph^2}(A_*)_{ab}(t)\right)\right|_{t=0}=0,
\quad
\left.\frac{d^2}{dt^2}\left(\frac{z^2}{\ph^2}(A_*)_{ab}(t)\right)\right|_{t=0}
=(A_*)''_{ab}(0)
+\left.\left(\frac{z^2}{\ph^2}\right)''\right|_{t=0} A_{ab}.
\]

In local coordinates $x(t)=(x^1(t),\ldots,x^n(t))$, 
\begin{align*}
\left.\frac{d^2}{dt^2}
\widetilde F\!\left(\bar A_*(t),Y,x(t)\right)\right|_{t=0}
&=\sum_{a,b,c,d}\widetilde F^{ab,cd}
 \left.\left(\frac{z^2}{\ph^2}(A_*)_{ab}\right)'\right|_{t=0}
 \left.\left(\frac{z^2}{\ph^2}(A_*)_{cd}\right)'\right|_{t=0} \\
&\quad+2\sum_{a,b,j}\widetilde F^{ab,x^j}
 \left.\left(\frac{z^2}{\ph^2}(A_*)_{ab}\right)'\right|_{t=0}
 \dot x^j(0) \\
&\quad+\sum_{a,b}\widetilde F^{ab}
 \left.\left(\frac{z^2}{\ph^2}(A_*)_{ab}\right)''\right|_{t=0}
 +\sum_{j,k}\widetilde F_{x^jx^k}\dot x^j(0)\dot x^k(0)
 +\sum_j\widetilde F_{x^j}\ddot x^j(0).
\end{align*}
All derivatives of $\widetilde F$ are evaluated at the common
initial matrix $A$, fixed $Y$, and $x_1$. 
Then
\begin{equation}\label{eq:claim-second}
\widetilde F^{x_i,x_i}-\widetilde F_D^{x_i,x_i}
=\sum_{a,b}\widetilde F^{ab}
\bigl((A_{\bn})''_{ab}(0)-(A_D)''_{ab}(0)\bigr).
\end{equation}
Moreover, applying the chain rule to \(\widetilde F^{ab}\), we obtain $\widetilde F^{ab,x_i}=\widetilde F_D^{ab,x_i}$.
Differentiating $(A_D)'_{ab}$ once more, and using $\bn-D=0$ at $x_1$, gives
\begin{align*}
&(A_{\bn})''_{ab}(0)-(A_D)''_{ab}(0)\\
&\quad=A\bigl([\bn_{\be_i}(\bn-D)](\be_i,\be_a),\be_b\bigr)
+A\bigl(\be_a,[\bn_{\be_i}(\bn-D)](\be_i,\be_b)\bigr).
\end{align*}
Since $A=\widetilde U^\eps(x_1)$ and $\widetilde F^{ab}$ are
locally uniformly bounded, \eqref{eq:claim-second} and
\eqref{eq:connestimate} imply
\[
\abs{\widetilde F^{x_i,x_i}-\widetilde F_D^{x_i,x_i}}
\le C\abs{\bn_{\be_i}(\bn-D)}\le C\widetilde U_{ii}.
\]
This proves both assertions.
\end{proof}

By \eqref{eq:claim}, the last two terms in \eqref{eq:decomposition} are bounded by $C\widetilde{U}_{ii}$. In the last term of \eqref{eq:big-expansion}, retain explicitly the terms whose matrix indices $a,b,c,d$ belong to $G$. We then obtain 
\begin{equation}\label{eq:G-bracket}
\begin{aligned}
&\sum_i\phi^{ii}\Bigl\{
 \widetilde{F}^{ab,cd}
 \widetilde{U}_{ab;i}\widetilde{U}_{cd;i}
 +2\widetilde{F}^{ab,x_i}\widetilde{U}_{ab;i}
 +\widetilde{F}^{x_i,x_i}
 \Bigr\} \\
&=\sum_{i\in B}\Bigl(
 \sg_l(G)
 +\frac{\sg_1(B\mid i)^2-\sg_2(B\mid i)}
       {\sg_1(B)^2}
 \Bigr) \\
&\quad\times\Biggl\{
 \sum_{a,b,c,d\in G}
 \widetilde{F}^{ab,cd}
 \widetilde{U}_{ab;i}\widetilde{U}_{cd;i}
 +2\sum_{a,b\in G}
 \widetilde{F}_D^{ab,x_i}\widetilde{U}_{ab;i}
 +\widetilde{F}_D^{x_i,x_i}
 \Biggr\} \\
&\quad
 +O\Bigl(
 \phi_\eps
 +\sum_{i,j\in B}\abs{\bn\widetilde{U}_{ij}}
 \Bigr).
\end{aligned}
\end{equation}
The derivatives in this bracket are taken with $Y$ fixed, exactly as required by condition \eqref{eq:struct}.

We now return to the first term of \eqref{eq:big-expansion}. Note that $\phi^{ij,km}=\sg_{l+1}^{ij,km}+q_{\eps}^{ij,km}$. 
For $i\in B$, $j\in G$, and every derivative direction $\alpha$, 
\begin{equation}\label{eq:spatial-codazzi}
\widetilde{U}_{ij;\alpha}=\widetilde{U}_{i\alpha;j}=\widetilde{U}_{\alpha j;i}.
\end{equation}
At $x_1$,
\begin{equation}\label{eq:V}
V_{i\alpha}=\sigma_{1}(B)\widetilde{U}_{ii;\alpha}
-(\widetilde{U}_{ii}+\eps)\sum_{j\in B}\widetilde{U}_{jj;\alpha},\qquad i\in B.
\end{equation}
Since $\sigma_{l-1}((\widetilde{U}^\eps)\mid ij)=O(\phi_\eps)$ for
$i,j\in G$, $i\ne j$, for $\alpha,\beta$ fixed, 
\[
\begin{aligned}
\sum_{i,j,k,m} \sigma_{l+1}^{ij,km}
\widetilde{U}_{ij;\alpha}\widetilde{U}_{km;\beta}
&=\sum_{i\ne k}\sigma_{l+1}^{ii,kk}
 \widetilde{U}_{ii;\alpha}\widetilde{U}_{kk;\beta}
 +\sum_{i\ne j}\sigma_{l+1}^{ij,ji}
 \widetilde{U}_{ij;\alpha}\widetilde{U}_{ji;\beta} \\
&=\sum_{i\ne k}\sigma_{l-1}(\widetilde{U}^\eps\mid ik)
 \widetilde{U}_{ii;\alpha}\widetilde{U}_{kk;\beta}
 -\sum_{i\ne j}\sigma_{l-1}(\widetilde{U}^\eps\mid ij)
 \widetilde{U}_{ij;\alpha}\widetilde{U}_{ji;\beta} \\
&=O\left(
 \phi_\eps+\sum_{i,j\in B}\lvert\bn\widetilde{U}_{ij}\rvert
 \right)
 -2\sum_{i\in B,j\in G}
 \sigma_{l-1}(G\mid j)
 \widetilde{U}_{ij;\alpha}\widetilde{U}_{ij;\beta}.
\end{aligned}
\]

As $\sigma_{l-1}(G\mid j)
=\dfrac{\sigma_l(G)}{\widetilde{U}_{jj}+\eps}$,
$j\in G$, we have
\[
\sigma_{l+1}^{ij,km}
\widetilde{U}_{ij;\alpha}\widetilde{U}_{km;\beta}
=O\left(
 \phi_\eps+\sum_{i,j\in B}\lvert\bn\widetilde{U}_{ij}\rvert
 \right)
 -2\sigma_l(G)\sum_{i\in B,j\in G}
 \frac{1}{\widetilde{U}_{jj}+\eps}
 \widetilde{U}_{ij;\alpha}\widetilde{U}_{ij;\beta}.
\]

By Lemma 2.5,
\[
\begin{aligned}
\sum_{i,j,k,m}q_{\eps}^{ij,km}
\widetilde{U}_{ij;\alpha}\widetilde{U}_{km;\beta}
&=O\left(
 \phi_\eps+\sum_{i,j\in B}\lvert\bn\widetilde{U}_{ij}\rvert
 \right)
 -2\sum_{i\in B,j\in G}
 \frac{\sigma_1^2(B\mid i)-\sigma_2(B\mid i)}
      {\sigma_1^2(B)(\widetilde{U}_{jj}+\eps)}
 \widetilde{U}_{ij;\alpha}\widetilde{U}_{ij;\beta} \\
&\quad
 -\frac{\sum_{i\in B}V_{i\alpha}V_{i\beta}}
       {\sigma_1^3(B)}
 -\frac{1}{\sigma_1(B)}
 \sum_{\substack{i,j\in B\\ i\ne j}}
 \widetilde{U}_{ij;\alpha}\widetilde{U}_{ji;\beta}.
\end{aligned}
\]

This gives
\begin{equation}\label{eq:new-formula}
\begin{aligned}
&\sum \widetilde{F}^{\alpha\beta}\phi^{ij,km}
 \widetilde{U}_{ij;\alpha}\widetilde{U}_{km;\beta}
 =O\left(
 \phi_\eps+\sum_{i,j\in B}\lvert\bn\widetilde{U}_{ij}\rvert
 \right) \\
&\quad {}-\sum_{\alpha,\beta}\widetilde{F}^{\alpha\beta}
 \Biggl[
 \frac{\sum_{i\in B}V_{i\alpha}V_{i\beta}}{\sigma_1^3(B)}
 +\frac{\sum_{\substack{i,j\in B\\ i\ne j}}
         \widetilde{U}_{ij;\alpha}\widetilde{U}_{ji;\beta}}
        {\sigma_1(B)} \\
&\qquad {}+2\sum_{\substack{i\in B\\ j\in G}}
 \left(
 \sigma_l(G)
 +\frac{\sigma_1^2(B\mid i)-\sigma_2(B\mid i)}{\sigma_1^2(B)}
 \right)
 \frac{\widetilde{U}_{ij;\alpha}\widetilde{U}_{ji;\beta}}
      {\widetilde{U}_{jj}+\eps}
 \Biggr].
\end{aligned}
\end{equation}
Combining \eqref{eq:X-error}, \eqref{eq:G-bracket},
and \eqref{eq:new-formula}, one reduces
\eqref{eq:chain-contraction} to
\begin{equation}\label{eq:combined}
\begin{aligned}
\widetilde{F}^{\alpha\beta}(\phi_\eps)_{;\alpha\beta}
&\le C\Bigl(
 \phi_\eps+\sum_{i,j\in B}\abs{\bn\widetilde{U}_{ij}}
 \Bigr)
 -\sum_{\alpha,\beta}\widetilde{F}^{\alpha\beta}
 \Biggl[
 \frac{\sum_{i\in B}V_{i\alpha}V_{i\beta}}
      {\sg_1(B)^3}
 +\frac{\sum_{\substack{i,j\in B\\ i\ne j}}
        \widetilde{U}_{ij;\alpha}\widetilde{U}_{ij;\beta}}
       {\sg_1(B)}
 \Biggr] \\[3pt]
&{}-\sum_{i\in B}\left(
 \sg_l(G)
 +\frac{\sg_1(B\mid i)^2-\sg_2(B\mid i)}
       {\sg_1(B)^2}
 \right)
 \Biggl\{
 \sum_{a,b,c,d\in G}
 \widetilde{F}^{ab,cd}
 \widetilde{U}_{ab;i}\widetilde{U}_{cd;i} \\[3pt]
&\qquad
 +2\sum_{a,c,j\in G}\widetilde{F}^{ac}
 \frac{\widetilde{U}_{aj;i}\widetilde{U}_{cj;i}}
      {\widetilde{U}_{jj}+\eps}
 +2\sum_{a,b\in G}
 \widetilde{F}_D^{ab,x_i}\widetilde{U}_{ab;i}
 +\widetilde{F}_D^{x_i,x_i}
 \Biggr\} \\[3pt]
&{}-\sum_i\phi^{ii}\sum_{a\ne i}
 \widetilde{F}^{aa}\Rc_{aiai}
 (\widetilde{U}_{aa}-\widetilde{U}_{ii}).
\end{aligned}
\end{equation}
We now apply condition \eqref{eq:struct} to the expression in braces in \eqref{eq:combined}. Fix $i\in B$ and let $x(t)$ be the $g_Y$-geodesic satisfying $x(0)=x_1$ and $\dot x(0)=\be_i$. Use a $D$-parallel $g_Y$-orthonormal frame along this geodesic, and define the positive matrix $A(t)$ by
\[
\begin{aligned}
A(t)={}&\bigl(\widetilde{U}_{ab}+\eps\delta_{ab}\bigr)_{a,b\in G}^{-1}\\
&-t\bigl(\widetilde{U}_{ab}+\eps\delta_{ab}\bigr)_{a,b\in G}^{-1}
(\widetilde{U}_{ab;i})_{a,b\in G}
\bigl(\widetilde{U}_{ab}+\eps\delta_{ab}\bigr)_{a,b\in G}^{-1}.
\end{aligned}
\]
Differentiating $A(t)A(t)^{-1}=I$ at $t=0$ gives
\begin{equation}\label{eq:inverse-curve-derivatives}
\begin{aligned}
A(0)^{-1}&=(\widetilde{U}_{ab}+\eps\delta_{ab})_{a,b\in G},\\
\left.\frac{d}{dt}\bigl(A(t)^{-1}\bigr)_{ab}\right|_{t=0}&=\widetilde{U}_{ab;i},\\
\left.\frac{d^2}{dt^2}\bigl(A(t)^{-1}\bigr)_{ab}\right|_{t=0}
&=2\sum_{c,j\in G}\widetilde{U}_{ac;i}
\left[(\widetilde{U}_{rs}+\eps\delta_{rs})_{r,s\in G}^{-1}\right]^{cj}\widetilde{U}_{jb;i}.
\end{aligned}
\end{equation}
Because $\widetilde{U}$ is diagonal at the point,
\begin{equation}\label{eq:inverse-curve-diagonal}
\left[(\widetilde{U}_{rs}+\eps\delta_{rs})_{r,s\in G}^{-1}\right]^{cj}
=\frac{\delta_{cj}}{\widetilde{U}_{jj}+\eps},
\end{equation}
and hence
\begin{equation}\label{eq:inverse-curve-contract}
\widetilde{F}^{ab}\left.\frac{d^2}{dt^2}\bigl(A(t)^{-1}\bigr)_{ab}\right|_{t=0}
=2\sum_{a,c,j\in G}\widetilde{F}^{ac}\frac{\widetilde{U}_{aj;i}\widetilde{U}_{cj;i}}{\widetilde{U}_{jj}+\eps}.
\end{equation}
When applying (1.4), we use $A(t)^{-1}$ as the $G\times G$ block
and keep all other entries fixed at their values in
$\widetilde U^\eps(x_1)$. The curve $A(t)$ is affine and $x(t)$ is a $g_Y$-geodesic. Hence condition \eqref{eq:struct}, with $Y$ fixed, gives
\begin{align*}
0
&\le \left.\frac{d^2}{dt^2}
 \widetilde{F}(A(t)^{-1},Y,x(t))\right|_{t=0}
 = \sum_{a,b,c,d\in G}\widetilde{F}^{ab,cd}
 \left.\frac{d(A(t)^{-1})_{ab}}{dt}
 \frac{d(A(t)^{-1})_{cd}}{dt}\right|_{t=0} \\
&\quad + \sum_{a,b\in G}\widetilde{F}^{ab}
 \left.\frac{d^2(A(t)^{-1})_{ab}}{dt^2}\right|_{t=0}
 + 2\sum_{a,b\in G}\widetilde{F}_D^{ab,x_i}
 \left.\frac{d(A(t)^{-1})_{ab}}{dt}\right|_{t=0}
 + \widetilde{F}_D^{x_i,x_i}.
\end{align*}
Substituting \eqref{eq:inverse-curve-derivatives}, \eqref{eq:inverse-curve-diagonal}, and \eqref{eq:inverse-curve-contract} into the preceding convexity inequality gives
\begin{equation}\label{eq:convexity-bracket}
\begin{aligned}
0\le{}&\sum_{a,b,c,d\in G}\widetilde{F}^{ab,cd}\widetilde{U}_{ab;i}\widetilde{U}_{cd;i}
+2\sum_{a,c,j\in G}\widetilde{F}^{ac}\frac{\widetilde{U}_{aj;i}\widetilde{U}_{cj;i}}{\widetilde{U}_{jj}+\eps}\\
&+2\sum_{a,b\in G}\widetilde{F}_D^{ab,x_i}\widetilde{U}_{ab;i}+\widetilde{F}_D^{x_i,x_i}.
\end{aligned}
\end{equation}
By uniform ellipticity, there exists a constant
$\delta_0>0$, such that
\begin{equation}\label{eq:next-formula}
\
(\widetilde{F}^{\alpha\beta})\ge\delta_0 I,
\qquad \forall x\in\mathcal{O}.
\end{equation}
Using \eqref{eq:convexity-bracket} and \eqref{eq:next-formula} in \eqref{eq:combined}, we then obtain
\[
\begin{aligned}
\sum_{\alpha,\beta}\widetilde{F}^{\alpha\beta}(\phi_\eps)_{;\alpha\beta}
&\le C\left(
 \phi_\eps+\sum_{i,j\in B}\lvert\bn\widetilde{U}_{ij}\rvert
 \right)
 -\delta_0\sum_{\alpha=1}^{n}
 \left[
 \frac{\sum_{i\in B}V_{i\alpha}^2}{\sigma_1^3(B)}
 +\frac{\sum_{\substack{i,j\in B\\ i\ne j}}
        \lvert\widetilde{U}_{ij;\alpha}\rvert^2}{\sigma_1(B)}
 \right] \\
&{}-\sum_i\phi^{ii}\sum_{a\ne i}
 \widetilde{F}^{aa}\Rc_{aiai}
 (\widetilde{U}_{aa}-\widetilde{U}_{ii}).
\end{aligned}
\]
By Lemma 3.3 in [3], for each $M\ge 1$, for any
$M\ge \gamma_i\ge \frac{1}{M}$, there is a constant $C$
depending only on $n$ and $M$ such that, $\forall\alpha$,
\[
\sum_{i,j\in B}\lvert\widetilde{U}_{ij;\alpha}\rvert
\le C\left(1+\frac{1}{\delta_0^2}\right)
 \left(
 \sigma_1(B)
 +\left|\sum_{i\in B}\gamma_i\widetilde{U}_{ii;\alpha}\right|
 \right)
 +\frac{\delta_0}{2}\left[
 \frac{\sum_{i\ne j\in B}\lvert\widetilde{U}_{ij;\alpha}\rvert^2}
      {\sigma_1(B)}
 +\frac{\sum_{i\in B}V_{i\alpha}^2}{\sigma_1^3(B)}
 \right].
\]

Set
\[
\gamma_j=\sigma_l(G)
 +\frac{\sigma_1^2(B\mid j)-\sigma_2(B\mid j)}{\sigma_1^2(B)},
 \qquad \forall j\in B.
\]
The Newton--Maclaurin inequality implies
\[
\sigma_l(G)+1\ge\gamma_j\ge\sigma_l(G),
\qquad \forall j\in B.
\]
Then there are positive constants $C_1$, $C_2$, independent of $\eps$, such that 
\begin{equation}\label{eq:phi-eps-estimate}
\sum_{\alpha,\beta} \widetilde{F}^{\alpha\beta}(\phi_\eps)_{;\alpha\beta}
\le C_1\bigl(\phi_\eps+\lvert\bn\phi_\eps\rvert\bigr)
-C_2\sum_{i,j\in B}\lvert\bn\widetilde{U}_{ij}\rvert
-\sum_i\phi^{ii}\sum_{a\ne i}
 \widetilde{F}^{aa}\Rc_{aiai}
 (\widetilde{U}_{aa}-\widetilde{U}_{ii}).
\end{equation}
By \eqref{eq:chain1} and \eqref{eq:phi-coeff}, 
\[
\sum_{i\in B}\gamma_i\widetilde U_{ii;\alpha}
=(\phi_\varepsilon)_{;\alpha}+O(\phi_\varepsilon).
\]
We now determine the sign of the curvature term before passing to the limit. 
For $i\in B$, every summand containing $\widetilde{U}_{ii}$ is $O(\phi_\eps)$ by \eqref{eq:bad-small}. For $i\in G$, $\phi^{ii}=O(\phi_\eps)$. Hence
\begin{equation}\label{eq:curvature-sign}
\begin{aligned}
&-\sum_i\phi^{ii}\sum_{a\ne i}\widetilde{F}^{aa}\Rc_{aiai}(\widetilde{U}_{aa}-\widetilde{U}_{ii})\\
&=-\sum_{i\in B}\left(\sg_l(G)
+\frac{\sg_1(B\mid i)^2-\sg_2(B\mid i)}{
\sg_1(B)^2}\right)
\sum_{a\in G}\widetilde{F}^{aa}\Rc_{aiai}\widetilde{U}_{aa}+O(\phi_\eps).
\end{aligned}
\end{equation}
Since
\[
\widetilde{F}^{aa}>0,\qquad\Rc_{aiai}\ge1,\qquad \widetilde{U}_{aa}>0\quad(i\in B,a\in G),
\]
its displayed principal part is nonpositive. If $B\ne\varnothing$ and $G\ne\varnothing$, this principal part is strictly negative after the eigenvalues indexed by $B$ are known to vanish identically. Combining \eqref{eq:phi-eps-estimate} and \eqref{eq:curvature-sign} gives
\begin{equation}\label{eq:curvature-final}
\begin{aligned}
\widetilde{F}^{\alpha\beta}(\phi_\eps)_{;\alpha\beta}
\le{}&C_1\bigl(\phi_\eps+\abs{\bn\phi_\eps}\bigr)-C_2\sum_{i,j\in B}|\bn\widetilde{U}_{ij}|\\\
&-\sum_{i\in B}\left(\sg_l(G)
+\frac{\sg_1(B\mid i)^2-\sg_2(B\mid i)}{
\sg_1(B)^2}\right)
\sum_{a\in G}\widetilde{F}^{aa}\Rc_{aiai}\widetilde{U}_{aa}.
\end{aligned}
\end{equation}
After dropping the last two nonpositive terms, this is the required differential inequality.
Taking $\eps \to 0$, we obtain 
\begin{equation}\label{eq:viscosity}
\widetilde{F}^{\alpha\beta}\phi_{;\alpha\beta}
\le C(\phi+\abs{\bn\phi}),\qquad\phi\ge0,\qquad\phi(x_0)=0.
\end{equation}
The strong minimum principle, applied to the nonnegative function $\phi$  gives $\phi\equiv0$  in $\mathcal O$.  And
\[
\widetilde{U}\ge0,\quad\phi=0\quad\Longrightarrow\quad\sg_{l+1}(\widetilde{U})=0
\quad\Longrightarrow\quad\operatorname{rank}\widetilde{U}\le l.
\]
By the definition of $l$, the reverse inequality holds. Hence
\[
S=\{x\in\Om:\operatorname{rank}\widetilde{U}(x)=l\}
=\{x\in\Om:\sg_{l+1}(\widetilde{U}(x))=0\}
\]
satisfies
\[
S\ne\varnothing,\qquad S\text{ is open},\qquad S\text{ is closed}.
\]
The connectedness of $\Om$ gives $S=\Om$. Since $\widetilde{U}=\ph^{-1}U[\ph]$, the tensor $U[\ph]$ has constant rank.

Suppose this rank were $l$ with $0<l<n$. Once $\phi\equiv0$, return to \eqref{eq:curvature-final} and retain, rather than discard, its curvature term. Letting $\eps \to 0$ in that inequality gives
\begin{equation}\label{eq:strict-curv}
\widetilde{F}^{\alpha\beta}\phi_{;\alpha\beta}
\le C(\phi+\abs{\bn\phi})
-\sum_{i\in B}\left(\sg_l(G)+\frac{\abs{B}-1}{2\abs{B}}\right)
\sum_{a\in G}\widetilde{F}^{aa}\Rc_{aiai}\widetilde{U}_{aa}.
\end{equation}
Indeed, $\widetilde{U}_{rr}=0$ for every $r\in B$, so $\widetilde{U}_{rr}+\eps=\eps$ for $r\in B$ and
\[
\frac{\sg_1(B\mid i)^2-\sg_2(B\mid i)}{
\bigl(\sum_{r\in B}(\widetilde{U}_{rr}+\eps)\bigr)^2}
=\frac{\abs{B}-1}{2\abs{B}}.
\]
Since $0<l<n$, both $B$ and $G$ are nonempty. On a sufficiently small neighborhood,
\[
\min_{i\in B}\left(\sg_l(G)+\frac{\abs{B}-1}{2\abs{B}}\right)>0,\quad
\min_{a\in G}\widetilde{F}^{aa}>0,\quad \Rc_{aiai}\ge1,\quad\min_{a\in G}\widetilde{U}_{aa}>0.
\]
Therefore
\[
\sum_{i\in B}\left(\sg_l(G)+\frac{\abs{B}-1}{2\abs{B}}\right)
\sum_{a\in G}\widetilde{F}^{aa}\Rc_{aiai}\widetilde{U}_{aa}>0.
\]
Because $\phi=0$, $\bn\phi=0$, and $\bn^2\phi=0$, substitution into \eqref{eq:strict-curv} gives
\[
0\le-\sum_{i\in B}\left(\sg_l(G)+\frac{\abs{B}-1}{2\abs{B}}\right)
\sum_{a\in G}\widetilde{F}^{aa}\Rc_{aiai}\widetilde{U}_{aa}<0,
\]
which is impossible. Thus $l\in\{0,n\}$, or equivalently,
\[
U[\ph]\equiv0\qquad\text{or}\qquad U[\ph]>0.
\]
\end{proof}

\begin{proof}[Proof of Theorem~\ref{thm:main}]
The local constant-rank conclusion and the alternative in Theorem~\ref{thm:differential} apply to $\widetilde{U}=\ph^{-1}U[\ph]$. Since $U[\ph]=\ph\widetilde{U}$ and $\ph>0$, these two tensors have
the same rank and nullspace. Therefore, the conclusions for
$\widetilde{U}$ also hold for $U[\ph]$.
\end{proof}

\section{An example: the $L_p$-Christoffel--Minkowski equation}

We illustrate the application of Theorem~\ref{thm:main} by considering the
$L_p$-Christoffel--Minkowski equation
\begin{equation}\label{eq:Lp-CM}
\sigma_k(U[\varphi])=\varphi^{p-k}f(x),
\qquad U[\varphi]\ge0,
\end{equation}
where $1\le k\le n$, $p\in\mathbb R$, and
$f\in C^\infty(\Omega)$ is positive. We consider the equation in the \(k\)-admissible cone.

Since $\varphi>0$ and $f>0$, equation \eqref{eq:Lp-CM} is equivalent to
\begin{equation}\label{eq:Lp-root}
-\sigma_k(U[\varphi])^{-1/k}
+\varphi^{1-p/k}f(x)^{-1/k}=0.
\end{equation}
To fit the transformed operator defined in \eqref{eq:tildeF}, we use the
following positive normalization of the root equation:
\begin{equation}\label{eq:Lp-F}
F(r,\xi,z,x)
:=z^{-1}\left(-\sigma_k(r)^{-1/k}
+z^{1-p/k}f(x)^{-1/k}\right),
\end{equation}
where $r$ is the symmetric-matrix variable and $\xi$ denotes the gradient
variable. The operator is independent of $\xi$. This normalization does not
change the solution set, since $z>0$.

On the cone $\Gamma_k$, the matrix derivative is
\[
F^{ij}
=\frac{1}{kz}\sigma_k(r)^{-1-1/k}\sigma_k^{ij}(r)>0.
\]
Thus $F$ is strictly elliptic in the $k$-admissible cone.

Fix $Y\in\mathbb H^{n+1}$ and put
\[
z(x)=-\langle Y,(x,1)\rangle>0,
\qquad \xi=\nabla z,
\qquad g_Y=z^{-2}g.
\]
By \eqref{eq:tildeF},
\begin{align}
\widetilde F(A,Y,x)
&=z^{-1}\left(-\sigma_k(z^{-1}A)^{-1/k}
+z^{1-p/k}f(x)^{-1/k}\right)\nonumber\\
&=-\sigma_k(A)^{-1/k}+z^{-p/k}f(x)^{-1/k}.
\label{eq:Lp-transformed}
\end{align}
Consequently,
\begin{equation}\label{eq:Lp-transformed-inverse}
\widetilde F(A^{-1},Y,x)
=-\sigma_k(A^{-1})^{-1/k}
+z^{-p/k}f(x)^{-1/k}.
\end{equation}

We now verify the structural condition \eqref{eq:struct}. For a positive
definite matrix $A$, the elementary symmetric-function identity gives
\[
\sigma_k(A^{-1})
=\frac{\sigma_{n-k}(A)}{\sigma_n(A)},
\]
where $\sigma_0=1$. Hence
\begin{equation}\label{eq:Lp-Newton-quotient}
\sigma_k(A^{-1})^{-1/k}
=\left(\frac{\sigma_n(A)}{\sigma_{n-k}(A)}\right)^{1/k}.
\end{equation}
The right-hand side is concave on the cone of positive definite matrices
(see \cite[Section~2.2]{HLX}). Hence
\(A\mapsto -\sigma_k(A^{-1})^{-1/k}\) is convex on this cone. 

The second term in \eqref{eq:Lp-transformed-inverse} depends only on $x$.
Let $D$ denote the Levi-Civita connection of $g_Y=z^{-2}g$, as in
Lemma~\ref{lem:freeze}. Since there are no mixed $(A,x)$ derivatives,
the joint convexity in $(A,x)$ reduces precisely to the convexity in $x$
of~$z^{-p/k}f(x)^{-1/k}$.

Thus a sufficient condition for \eqref{eq:struct} is
\begin{equation}\label{eq:Lp-convexity}
D^2\left(z^{-p/k}f^{-1/k}\right)\ge0
\qquad\text{for every fixed }Y\in\mathbb H^{n+1}.
\end{equation}
The left-hand side can be expanded as
\begin{align*}
0\le{}&D^2(f^{-1/k})
-\frac pk\,D\log z\otimes D(f^{-1/k})
-\frac pk\,D(f^{-1/k})\otimes D\log z\\
&+f^{-1/k}\left[
-\frac pk D^2\log z
+\frac{p^2}{k^2}D\log z\otimes D\log z
\right].
\end{align*}
Here all covariant derivatives and tensor products are taken with respect to
$g_Y$. The factor $z^{-p/k}$ has been omitted from the displayed inequality
because it is positive.

Whenever \eqref{eq:Lp-convexity} holds for every fixed $Y$, the structural
condition of Theorem~\ref{thm:main} is satisfied. Applying that theorem to
\eqref{eq:Lp-CM}, we conclude that $U[\varphi]$ has constant rank. Moreover,
\[
\sigma_k(U[\varphi])=\varphi^{p-k}f(x)>0
\]
excludes the possibility $U[\varphi]\equiv0$. The alternative in
Theorem~\ref{thm:main} therefore yields
\[
U[\varphi]>0\quad\text{throughout }\Omega.
\]
This completes the application of the constant-rank theorem to the
$L_p$-Christoffel--Minkowski equation.

\bigskip

\noindent
School of Mathematics, Hunan University,
Changsha 410082, P. R. China.

\medskip

\noindent
\textit{Email address:} \texttt{xulu@hnu.edu.cn}

\medskip

\noindent
\textit{Email address:} \texttt{zhouyutian@hnu.edu.cn}

\medskip

\noindent
\textit{Email address:} \texttt{ttliu@hnu.edu.cn}

\end{document}